\documentclass[12pt]{hjm1}

\usepackage{amsmath,amssymb,amscd,amsthm,amsxtra}
\usepackage[colorlinks,citecolor=red,pagebackref,hypertexnames=false]{hyperref}
    
\numberwithin{equation}{section}

\newtheorem{theorem}[equation]{Theorem}
\newtheorem{lemma}[equation]{Lemma}
\newtheorem{proposition}[equation]{Proposition}

\def\ip#1,#2.{\langle #1,#2\rangle}
\def\Ip#1,#2.{\bigl\langle#1,#2\bigr\rangle}
\def\IP#1,#2.{\Bigl\langle#1,#2\Bigr\rangle}
\def\norm#1.#2.{\lVert#1\rVert_{#2}}
\def\Norm#1.#2.{\bigl\lVert#1\bigr\rVert_{#2}}
\def\NOrm#1.#2.{\Bigl\lVert#1\Bigr\rVert_{#2}}
\def\NORm#1.#2.{\biggl\lVert#1\biggr\rVert_{#2}}
\def\NORM#1.#2.{\Biggl\lVert#1\Biggr\rVert_{#2}}

\def\abs#1{\lvert#1\rvert}

\def\ABs#1{\Bigl\lvert#1\Bigr\rvert}

\def\ind#1{\operatorname 1_{#1}}

 \def\Enl #1,#2,{\operatorname{\rm Enl}_{#1}(#2)}
 \def\emb #1.#2.{\operatorname{emb}(#1,#2)}
 \def\Emb#1..{\operatorname{emb}(#1)}

 \def\Rho#1{R_{#1}}

 \def\sh#1{\operatorname{sh}(#1)}

\begin{document}

\title
[Hankel Operators  and Product $BMO$]
{Hankel Operators in Several Complex Variables  and Product $BMO $}

\author{Michael Lacey}

\address{Michael Lacey\\
School of Mathematics\\
Georgia Institute of Technology\\
Atlanta,  GA 30332 USA}

\email{lacey@math.gatech.edu}

\author{Erin Terwilleger}

\address{
Erin Terwilleger\\
Department of Mathematics, U-3009\\
University of Connecticut\\
Storrs, CT 06269 USA}

\email{terwilleger@math.uconn.edu}

\date{\today}



\begin{abstract}   $H^2 (\otimes_1^n\mathbb C_+)$ denotes the Hardy space of square integrable functions
analytic in each variable separately.  Let $P^{\ominus}$ be the natural projection of $L^2 (\otimes_1^n\mathbb C_+)$
onto $\overline {H^2 (\otimes_1^n\mathbb C_+)}$.  A Hankel operator with symbol $b$ is the linear operator from
$H^2 (\otimes_1^n\mathbb C_+)$ to $\overline {H^2 (\otimes_1^n\mathbb C_+)}$ given by $\operatorname H_b \varphi=P^{\ominus}\overline b \varphi$.
We show that
\begin{equation*}
\norm \operatorname H_b ..\simeq \norm P^{\oplus}b.BMO (\otimes_1^n\mathbb C_+).,
\end{equation*}
where the right hand norm is S.-Y.~Chang and R.~Fefferman product $BMO$.
This fact has well known equivalences in terms of commutators and the weak factorization of $H^1 (\otimes_1^n\mathbb C_+)$.
In the case of two complex variables, this is due to  Ferguson and Lacey \cite{MR1961195}.
While the current proof is inductive, and one can take the one complex variable case as the basis step,
it is heavily influenced by the methods of Ferguson and Lacey.
 The induction is carried out with a particular form of a lemma due to Journ\'e
\cite{MR87g:42028}, which occurs implicitly in the work of J.~Pipher \cite{MR88a:42019}.
\end{abstract}

\maketitle

\section{Introduction} 

We characterize the boundedness of Hankel operators in three and more  complex variables in terms of 
the $BMO$ norm of the symbol of the operator.    In one complex variable, this and related facts
are  the circle of ideas around Nehari's theorem.
In the case of more than one complex variable, there are different types of Hankel operators, and we only consider the
so called little Hankel operators; little in the sense that the projection used in the definition is onto the smallest
natural choice of subspaces of $L^2(\mathbb R^n)$ to use.   The structure of these Hankel operators is
  more intricate  due to the more complicated structure
of the Hardy spaces $H^1$ in the product domain and their duals, as identified by S.-Y.~Chang and R.~Fefferman [\cite{chang},\cite{MR584078},
\cite{MR766959}].  Some of the tools that have proved to be so flexible and powerful in the one parameter situation apparently have no analog in the higher parameter case;  these spaces remain, to a significant degree, poorly understood.

We prove the natural statement about the boundedness of little Hankel operators in an arbitrary number of complex variables. Namely,
the Hankel operator with symbol $b$ is bounded iff the projection of $b$ into product Hardy space is in product $BMO$.
Central to this paper is the result of S.~Ferguson and M.~Lacey \cite{MR1961195} that established a
similar characterization for Hankel operators of two complex variables.    The current proof is inductive in nature, and one can
use the classical one variable statements of our theorem as the base case in the induction. 
In particular the methods of \cite{MR1961195} are not sufficient to prove the Theorem 
in this paper; the inductive argument is the  essential new argument in this paper.

Recall that  $L^2(\mathbb R)$ has the orthogonal decomposition $H^2(\mathbb R)\oplus
\overline {H^2(\mathbb R)}$.  Let $\operatorname P^\pm$ be the corresponding orthogonal projections onto the
analytic/antianalytic  spaces.

In $n$ variables, let $P_j^\pm$ be the same projections acting on the $j$th coordinate, $j\in\{1,2,\ldots,n\}$.
For functions $\sigma\,:\,\{1,2,\ldots,n\}\longrightarrow\{\pm\}$, let
\begin{equation*}
\operatorname P^\sigma=\prod_{j=1}^n \operatorname  P_j^{\sigma(j)}.
\end{equation*}
It is clear that  $L^2(\mathbb R^n)$ has the orthogonal decomposition into
\begin{equation*}
L^2(\mathbb R^n)=\oplus_{\sigma}\operatorname P^\sigma L^2(\mathbb R^n).
\end{equation*}
We take $\oplus$ to be  the function from $\{1,2,\ldots,n\}$ that is identically $+$.
It is clear that $\operatorname P^{\oplus}L^2(\mathbb R^n)=H^2 (\otimes_1^n\mathbb C_+)$.

For a function $b$, we set the {\em Hankel operator with symbol $b$} to be
$\operatorname H_b f=\operatorname P^{\ominus} \overline b  f$, defined as a map from $H^2 (\otimes_1^n\mathbb C_+)$ to $\overline {H^2 (\otimes_1^n\mathbb C_+)}$.
$\operatorname P^{\ominus}$ is the projection from $L^2(\mathbb R^n)$ onto $\overline {H^2 (\otimes_1^n\mathbb C_+)}$.
Clearly, this
operator depends only on  $\operatorname P^{\oplus} b$.

\begin{theorem}\label{t.hankel}
  We have the equivalence of norms
\begin{equation}\label{e.hankel}
\norm \operatorname H_b ..\simeq{} \norm \operatorname P^{\oplus} b.BMO (\otimes_1^n\mathbb C_+) . .
\end{equation}
\end{theorem}

Here $BMO (\otimes_1^n\mathbb C_+)$ is the (analytic) Bounded Mean Oscillation space, dual to $H^1 (\otimes_1^n\mathbb C_+)$, as identified by S.-Y.~Chang and R.~Fefferman.
This theorem has two well known equivalences.  One is in terms of the commutator
\begin{equation*}
C_b:=[\cdots[\,[\operatorname M_b,\operatorname H_1],\operatorname H_2],\cdots, \operatorname H_n],
\end{equation*}
in which $\operatorname M_b$ is the operator of pointwise multiplication by $b$, and $\operatorname H_j$ denotes the
Hilbert transform computed in the $j$th coordinate.  The commutator is a sum of $2^n$
Hankel operators, each coming from one of the $2^n$ orthants of $\mathbb R^n$.
In particular, if the \emph {signature} of $\sigma$ is $\text{sgn}(\sigma)=\prod_{j=1}^n \sigma(j)$,
a straightforward computation shows that
$
C_b={}  -2^n
\sum_\sigma \text{sgn}(\sigma) \operatorname P^{-\sigma}\operatorname M_b \operatorname P^\sigma.$
Thus, the upper bound for the Hankel operators $\norm \operatorname H_b..\lesssim\norm b.BMO.$ immediately
extends to an $L^2$ operator norm for the commutators.  Conversely, assuming the commutator is bounded on $L^2 (\otimes_1^n\mathbb C_+)$, a number
of Hankel operators
with the same symbol are  also bounded.  Namely, the Hankel operators are  from $\operatorname P^{\sigma}L^2 $ to $\operatorname P^{-\sigma }L^2$.
Thus the lower bound follows.  That is, we have
\begin{equation*}
\norm C_b .2.\simeq\norm b.BMO (\otimes_1^n\mathbb C_+) . .
\end{equation*}
The latter space is the {\em real} $BMO$ space.

A second equivalence is in essence a dual statement to the estimates above, and hence is a statement about $H^1 (\otimes_1^n\mathbb C_+)$.
It gives us a {\em weak factorization} result for that space, namely
\begin{equation} \label{e.weak}
H^1 (\otimes_1^n\mathbb C_+)=H^2 (\otimes_1^n\mathbb C_+)\widehat\odot {}H^2 (\otimes_1^n\mathbb C_+),
\end{equation}
where the right hand side is the projective tensor product of $H^2$.   This equality plays a role in our proof, and so we return to it below.

For the proof, the  strategy is one of induction on the number of parameters in a manner analogous to
the overall strategy of Ferguson and Lacey \cite{MR1961195}, which addresses the two parameter case.
The upper bound on $\norm \operatorname H_b..$ in  (\ref{e.hankel})
is in fact easy to obtain, a fact which
is easiest to see via the trivial inclusion in (\ref{e.weak}).
Thus, the real difficulty in our theorem lies in the lower bound
on $\norm \operatorname H_b..$.   Here, there is a bound which follows from the one parameter theory, namely that the
operator norm of $\operatorname H_b$ is bounded below by the ``rectangular $BMO$'' norm of $b$.
It is well known that the rectangular $BMO$ norm is
essentially smaller than the $BMO$ norm.  Ferguson and Lacey \cite{MR1961195} showed how to
use the Journ\'e Lemma \cite{MR87g:42028} to pass from this essentially smaller norm to the $BMO$ norm of Chang and Fefferman.

A formulation of the Journ\'e Lemma for rectangles in three and higher parameters is due to J.~Pipher \cite{MR88a:42019},
but the direct application of this lemma cannot succeed in a proof of our theorem.
The reasons are both   technical and heuristic.  Relying on just the rectangular $BMO$ norm in three and more parameters does not take
advantage of the subtle way that the $n$ parameter
$BMO$ space is built up from the $n - 1$ parameter space.

 We find that this point of view, and a form of Journ\'e's Lemma we need,
 as stated in Section~\ref{s.journe},  are  implicit in the paper of J.~Pipher.
To use the Journ\'e Lemma, we need to make a definition of $BMO_{-1} (\otimes_1^n\mathbb C_+)$, which is applied to a function in the $n$ parameter setting.
Our induction argument then, in proving the lower bound in the $n$ parameter setting, is to first derive the weaker bound of
$BMO_{-1} (\otimes_1^n\mathbb C_+)$.  And then prove the correct $BMO$ bound, assuming that the $BMO_{-1} (\otimes_1^n\mathbb C_+)$ norm of the symbol is sufficiently small.

By $A\lesssim{}B$ we mean that there is an absolute constant $K$ for which $A\le{}KB$.
$K$ is allowed to depend upon relevant parameters.

We are indebted to J.~Pipher for sharing some of her insights into the Journ\'e Lemma, 
and to the referee for a quick and helpful report.


\section{The Upper Bound}

The upper bound $\norm \operatorname H_b .2. {}\lesssim\norm b.BMO (\otimes_1^n\mathbb C_+).$ can be seen by a soft proof.
Consider the Hankel operator $\operatorname H_b$,
\begin{align}\nonumber
\norm \operatorname H_b..={}&
\sup_{\substack{f ,g\in  H^2 (\otimes_1^n\mathbb C_+) \\
\norm f.2.=1, \norm g.2.=1}}
\int  (\operatorname P^{\ominus} \overline b f)  g \; dx
\\{}={}& \sup_{\substack{f ,g\in  H^2 (\otimes_1^n\mathbb C_+) \\
\norm f.2.=1, \norm g.2.=1}} \int\overline{ \operatorname P^{\oplus}b } {fg} \;
dx .\label{bnorm1} \end{align}
 Since the product of $H^2$ functions is in
$H^1$, we see that the integral  above admits the upper bound of
$\norm \operatorname P^{\oplus}b.BMO (\otimes_1^n\mathbb C_+).$. This is the upper half of Theorem~\ref{t.hankel}.

We turn to the weak factorization result.
 For $A$, $B$ closed subspaces of $L^2(\mathbb R^n)$, we 
 define the projective tensor product  $A \widehat\odot  B \subseteq
 L^1(\mathbb R^n)$ by
 $$
 A \widehat\odot  B:= \{ h=\sum_{j=1}^{\infty} f_jg_j
 \,\,\Big| \,\,(f_j) \subseteq A, (g_j) \subseteq B, \textup{ and
 }\sum_{j=1}^{\infty}\norm f_j.2. \norm g_j.2. < \infty\Big\}.
 $$
 Observe that   that (\ref{bnorm1})  implies that the Hankel operator with symbol $b$ is bounded  if and only if the function $\operatorname P^{\oplus}b$ is in the dual of
 $H^2 (\otimes_1^n\mathbb C_+) \widehat\odot  H^2 (\otimes_1^n\mathbb C_+)$.   Therefore, the weak factorization equivalence  (\ref{e.weak}) is equivalent to our main theorem.   That is, we have the equivalence
 \begin{equation*}
 \norm \operatorname H_b..\simeq \norm b .(H^2 (\otimes_1^n\mathbb C_+) \widehat\odot  H^2 (\otimes_1^n\mathbb C_+))^*..
 \end{equation*}



\section{Wavelets, $BMO (\otimes_1^n\mathbb C_+)$ and $BMO_{-1} (\otimes_1^n\mathbb C_+)$}  

We begin with some preliminary definitions and calculations in the one parameter setting
which carry over naturally  to the higher parameter setting.
The proofs in the rest of the paper use analytic  wavelets constructed by Y. Meyer \cite{MR98e:42001} which are
compact in frequency.  Let $w$ be
 a Schwartz function with $\norm w.2.=1$ and $\widehat w (\xi)$ supported on $ [2/3,8/3]$.  Therefore the
 wavelets and projections  have the nice decay estimates
$$
|w(x)|{}\lesssim{} (1+|x|)^{-n} \quad\textup{ for } n \ge 1.
$$

Let $\mathcal D$ denote the dyadic intervals on $\mathbb R$.  For an interval $I \in \mathcal D$, define
$$
w_{I}(x):=|I|^{-\frac12} w\left(\frac{x-c(I)}{|I|}\right),
$$
where $c(I)$ denotes the center of $I$.
Note that the functions $w_I(x)$ are well localized to the interval $I$.  Indeed,
$$
|w_I(x)|  {}\lesssim{} |I|^{-\frac12}\Big
(\frac{1+\textup{dist}(x,I)}{|I|}\Big)^{-n} \textup{ for } n \ge 1.
$$
Y.~Meyer has shown that we can choose $w$ so that   $\{w_I\}_{I \in \mathcal D}$ form an orthonormal basis on $H^2 (\otimes_1^n\mathbb C_+)$.  Another useful property of these
functions is that we have the Littlewood-Paley inequalities,
$$
\NOrm {\sum_I \ip f,w_I.w_I} .p. \simeq \NOrm\Big( \sum_{I \in \mathcal D} \frac{|\ip f,w_I.|^2}{|I|}
1_I\Big)^{\!\!\frac12}.p.,
\qquad 1<p<\infty.
$$

We are now ready to define a characterization of product $BMO (\otimes_1^n\mathbb C_+)$  due to
 S.- Y. Chang and R. Fefferman \cite{MR584078}.
Let $\mathcal R=\mathcal  D^n=\otimes_{j=1}^n\mathcal D $ be the dyadic rectangles.
For a rectangle
 $R=\otimes_{j=1}^n R_j\in \mathcal R$, define
$$
v_R(x):=\prod_{j=1}^n w_{R_j}(x_j).
$$
 We say $f \in BMO (\otimes_1^n\mathbb C_+)$ iff
$$
\sup_{U} \bigg[|U|^{-1} \sum_{R \subset U} |\ip f,v_R.|^2
\bigg]^{\frac12} < \infty
$$
where $U$ is an open set in $\mathbb R^n$ of finite measure.\footnote{This is analytic $BMO$, as we are using
analytic wavelets.  Real $BMO$ has a similar definition, provided one uses wavelets that form a basis for $L^2(\mathbb R)$.}
We denote this supremum by $\norm f.BMO (\otimes_1^n\mathbb C_+).$.
It is a theorem of Chang and Fefferman that this definition coincides with the norm of the dual to $H^1 (\otimes_1^n\mathbb C_+)$.

We now define a weaker notion, which we dub $BMO_{-1} (\otimes_1^n\mathbb C_+)$.   For a collection of rectangles $\mathcal U\subset\mathcal R$,
set the {\em shadow } of $\mathcal U$, to be
\begin{equation*}
\operatorname{sh}(\mathcal U):=\bigcup_{R\in\mathcal U}R.
\end{equation*}
We say that $\mathcal U$ has $n - 1$ parameters iff there is a coordinate $1\le{}k\le{}n$ and a dyadic interval $I$ so that for all $R\in\mathcal U$, we have $R_k=I$. 
We then define
\begin{equation*} \label{e.n - 1}
\norm b.BMO_{-1} (\otimes_1^n\mathbb C_+).=\sup_{\text{$\mathcal U$, $n - 1$ parameters}}\Bigl[ \abs{\operatorname{sh}(\mathcal U)}^{-1} \sum_{R\in \mathcal U} \abs{\ip b,v_R.}^2 \Bigr]^{\frac12}.
\end{equation*}
Here, we note that the definition depends only upon the projection $\operatorname P^{\oplus}b$.
In two dimensions, this reduces to a notion that is just slightly weaker than the notion of rectangular $BMO$, which is well known to
be essentially smaller than the $BMO$ norm.


\section{The $BMO_{-1} (\otimes_1^n\mathbb C_+)$ lower bound}

An essential part of the argument is to use the induction
hypothesis to show that we have the lower bound \begin{equation*}
 \norm \operatorname H_b..\gtrsim{} \norm \operatorname P^{\oplus}
b.BMO_{-1} (\otimes_1^n\mathbb C_+).. \end{equation*} This amounts to the assertion that
\begin{equation}\label{e.n-1>} \norm b.(H^2 (\otimes_1^n\mathbb C_+) \widehat\odot 
H^2 (\otimes_1^n\mathbb C_+))^*.\gtrsim{}\norm b.BMO_{-1} (\otimes_1^n\mathbb C_+)., \end{equation} an inequality
we will demonstrate by relying on the truth of Theorem~\ref{t.hankel} in
the $n - 1$ parameter setting.

Given a symbol $b=b(x_1,x_2,\ldots,x_n)=b(x_1,x')$ of $n$ variables, we assume that $b$ is analytic in all variables and
has $\norm b.BMO_{-1} (\otimes_1^n\mathbb C_+).=1$.   We also take as given a set $\mathcal U$ of rectangles in $\mathcal D^n$ of
$n - 1 $ parameters.  Thus associated to $\mathcal U$ are the dyadic interval $I$ and the collection $\mathcal U^{(n - 1 )}\subset\mathcal D^{n - 1 }
$ as in the definition.   We assume that $\abs I=1$, $\abs{\sh {\mathcal U}}\simeq1$, and  for all $R\in\mathcal U$ we have
$\Rho1=I$ and $R=I\times R'$.

Our claim is then that the function \begin{equation*} \psi=\sum_{R\in\mathcal U}\ip b,
v_R.v_R \end{equation*} has $H^2 (\otimes_1^n\mathbb C_+) \widehat\odot  H^2 (\otimes_1^n\mathbb C_+)$ norm
${}\lesssim{}1$, which, together with $\ip \psi,b.=1$, certainly
proves (\ref{e.n-1>}). Thus it suffices to show the claim.

Now, since for each $R\in\mathcal U$, we have
$v_R(x_1,x')=w_{I}(x_1)v_{R'}(x')$, and \begin{equation*} \psi(x_1,x')={}
w_I(x_1)\sum_{R\in\mathcal U}\ip b,v_R. v_{R'}(x'):=w_I(x_1)\psi'(x'),
\end{equation*} we can utilize factorization results in both $x_1$ and $x'$.
For $x_1$, we use the classical inner outer factorization to
conclude that \begin{equation*} w_I=w_I^{(1)}w_I^{(2)},\qquad
\norm w_I^{(1)}.H^2({\mathbb C}_+).\norm
w_I^{(2)}.H^2({\mathbb C}_+).\lesssim{}1. \end{equation*} Concerning
the function $\psi'$, by our choice of $\mathcal U$ and the square function characterization
of the Hardy space, we observe that $\psi'$ has
$H^1(\otimes_{j=1}^{n-1}{\mathbb C}_+)$ norm at most a constant.
 By the induction hypothesis and, in particular,
the assertion that $H^1(\otimes_{j=1}^{n-1}{\mathbb C}_+)={}
H^2(\otimes_{j=1}^{n-1}{\mathbb C}_+)\widehat\odot 
H^2(\otimes_{j=1}^{n-1}{\mathbb C}_+)$, we can write \begin{equation*}
\psi'(x')=\sum_k \varphi_k(x')\phi_k(x'),\qquad \sum_k  \norm
\varphi_k.H^2(\otimes_{j=1}^{n-1}{\mathbb C}_+).\norm
\phi_k.H^2(\otimes_{j=1}^{n-1}{\mathbb C}_+).\lesssim1. \end{equation*} Hence,
writing \begin{equation*} \psi(x_1,x')=\sum_k
\bigl[w_I^{(1)}(x_1)\varphi_k(x')\bigr]\cdot\bigl[w_I^{(2)}
(x_1)\phi_k(x')\bigr] \end{equation*} we see that our claim holds. This
completes the proof of the lower bound.


\section{The $BMO$ lower bound}

An example of Carleson shows that the $BMO_{-1} (\otimes_1^n\mathbb C_+)$ bound is essentially smaller than the $BMO$ norm.
In particular, some tool is needed to pass to the larger norm.  That tool is a Journ\'e Lemma, which we detail in the next section.

We can assume that $\norm b.BMO (\otimes_1^n\mathbb C_+).=1$ and seek an absolute
lower bound on $\norm \operatorname H_b..$. In the course of the proof, we will
need absolute positive constants $\delta_{-1},\delta_{\text{journ\'e}},\delta_2$,
and $\delta_3$. Other parameters, termed ``diagonalization    parameters'',
are introduced to gain convergent geometric series.
The parameters used for these will be denoted with the letter $d$ with various subscripts.
We
 assume that $\norm b.BMO_{-1} (\otimes_1^n\mathbb C_+).<\delta_{-1}$, for
otherwise we have an absolute lower bound on $\norm \operatorname H_b..$.

Take a set of rectangles $\mathcal U$ which achieves the supremum in the
definition of the $BMO$ norm of $B$.  We can assume, after a harmless
dilation, that $\frac12<\abs{\operatorname{sh}(\mathcal U)}\le{}1$.   We
will show that \begin{equation} \label{e.lower} \norm \operatorname H_b
\overline {\operatorname P[\mathcal U]b}.2.\ge{}\delta_3. \end{equation} Here we use the notation
$\operatorname P[\mathcal U]=\sum_{R\in\mathcal U}v_R\otimes v_R$, and so
$\operatorname P[\mathcal U]b=\sum_{R\in\mathcal U}\ip b,v_R.v_R$. Establishing the lower
bound on the norm of the Hankel operator will require some careful
analysis which centers around a variety of paraproducts,
proper formulation, and application of a lemma due to Journ\'e which is specified in Section~\ref{s.journe}.

From the discussion in Section~\ref{s.journe}, there is a set
$V\supset\operatorname{sh}(\mathcal U)$, satisfying several conditions, among them 
$$
\abs
V\le{}(1+\delta_{\text{journ\'e}})\abs{\operatorname{sh}(\mathcal U)}\,.
$$
Take the collection of rectangles $\mathcal V$ and $\mathcal W$  to be \begin{gather*}
\mathcal V:=\{R\in\mathcal R\mid R\subset V,\ R\not\subset \sh{\mathcal U}\},
\\
\mathcal W:=\mathcal R-\mathcal U-\mathcal V.
\end{gather*}
We shall prove that for absolute $\delta_2>0$,
\begin{gather}   \label{e.u>}
\norm \operatorname P^{\ominus} \operatorname P[\mathcal U]b\overline {\operatorname P[\mathcal U]b}.2.\ge\delta_2,
\\ \label{e.v<}
\norm \operatorname P^{\ominus}\operatorname P[\mathcal V]b\overline {\operatorname P[\mathcal U]b}.2.\lesssim{}\delta_{\text{journ\'e}}^{1/2},
\\ \label{e.w<}
\norm \operatorname P^{\ominus}\operatorname P[\mathcal W]b\overline {\operatorname P[\mathcal U]b}.2.\le{}K_{\delta_{\text{journ\'e}}}\delta_{-1}.
\end{gather}
The last inequality holds with a constant that depends only on $\delta_{\text{journ\'e}}$.  Thus, fixing first $\delta_{\text{journ\'e}}$
sufficiently small and then $\delta_{-1}$ proves (\ref{e.lower}).

The first two estimates are trivial, as we indicate now.  First,
note that the Fourier transform of $\abs{ \operatorname P[\mathcal U] b}^2$ is
symmetric, so that \begin{equation*} \norm \operatorname P^{\ominus}
\operatorname P[\mathcal U]b\overline {\operatorname P[\mathcal U]b}.2.\ge2^{-n}\norm \operatorname P[\mathcal U]b.4.^2. \end{equation*} The
$L^4$ norm has a lower bound, due to the fact that we have taken
the shadow of $\mathcal U$ to have measure approximately one and the
validity of the Littlewood-Paley inequalities.  Thus, \begin{align*}
\tfrac14\le{}& \norm \operatorname P[\mathcal U]b .2.
\\{}={}& \sum_{R\in\mathcal U} \abs{\ip b,v_R.}^2
\\{}\le{}& \NOrm \Bigl[\sum_{R\in\mathcal U} \frac {\abs{\ip b,v_R.}^2}{\abs R } \ind R \Bigr]^{1/2}.4.
\\{}\lesssim{}&  \norm \operatorname P[\mathcal U]b .4..
\end{align*}
This proves (\ref{e.u>}).

Second, use   the control on the size of $V$ to see that \begin{equation*}
\norm \operatorname P[\mathcal U]b .2.^2+\norm \operatorname P[\mathcal V]b.2.^2=\abs {\sh{\mathcal  U}}+\norm
\operatorname P[\mathcal V]b.2.^2\le{}\abs V\le{}(1+\delta_{\text{journ\'e}})\abs
{\sh{\mathcal  U}} . \end{equation*} Thus, $\norm
\operatorname P[\mathcal V]b.2.^2\le{}\delta_{\text{journ\'e}}$. By the John-Nirenberg
inequality, we see that $\norm
\operatorname P[\mathcal V]b.4.\le{}\delta_{\text{journ\'e}}^{1/4}$.  Hence, we can prove
(\ref{e.v<}) as follows. \begin{equation*} \norm
\operatorname P^{\ominus}\operatorname P[\mathcal V]b\overline {\operatorname P[\mathcal U]b}.2.\le{}
    \norm \operatorname P[\mathcal U]b.4.\norm \operatorname P[\mathcal V]b.4.\lesssim{} \delta_{\text{journ\'e}}^{1/2}.
\end{equation*}

\subsection*{The Definitions of the Paraproducts}

The principal inequality is (\ref{e.w<}), and it requires a
sustained analysis to verify. It is imperative to observe that the
term $\operatorname H_{\operatorname P[\mathcal W]b}\operatorname P[\mathcal U]b$ has a sizable cancellation as a sum
over wavelets.  If $R$ and $R'$ are two dyadic rectangles with
$\abs{R_j}<8\abs{R_j'}$ for any $1\le{}j\le{}n$, then we would
have \begin{equation*} \operatorname P^{\ominus}v_{R'}\overline {v_{R}}=0. \end{equation*} This is due to the
fact that in the $j$th coordinate, the Fourier transform is not
supported in $\xi_j<0$. Thus, we can replace the definition of
$\mathcal W$ by: \begin{equation*} \mathcal W=\{R'\in \mathcal R\mid R'\not\subset V,\ \text{and
for some $R\in\mathcal U$, $\abs{R_j'}<8\abs{R_j}$ for all
$1\le{}j\le{}n$}\}. \end{equation*} It is also imperative to observe that
even with this restricted definition,  the shadows of $\mathcal U$ and
$\mathcal W$ will, in general, overlap.  This overlap will be controlled
by the Journ\'e Lemma and  additional orthogonality
considerations.

Nevertheless, the sum should be analyzed along the lines of a
product of two functions which are nearly supported on disjoint
sets.  The technique for doing this  is via sums known generically
as paraproducts.  In $n$ parameters, the paraproducts admit
different degeneracies, as measured in the amount of orthogonality
present in the sums.   It is the purpose of the following
definitions to quantify these paraproducts.

Given a subset $J\subset\{1,2,\ldots,n\}$, write $R'\prec_J{}R$ iff for indices $j\in J$, we have $8\abs{R'_j}<\abs {R_j}$,
 whereas for indices  $j\in \{1,2,\ldots,n\}-J$, we have $8^{-1}\abs{R'_j}\le\abs {R_j}\le{} 8\abs{R'_j}$.
Set
\begin{gather*}
\label{e.defX}
\mathcal X(J):=\{(R',R)\in \mathcal W\times\mathcal U\mid R'\prec_J R\}
\\ \label{e.defsigma}
\mathbb X(J):=\sum_{(R',R)\in\mathcal X(J)} \overline {\ip b,v_{R'}.v_{R'}}\ip b,v_{R}.v_{R}.
\end{gather*}
The remainder of the proof is devoted to the assertion that
\begin{equation} \label{e.ZX<}
\norm \mathbb X(J).2.\le{}K_{\delta_{\text{journ\'e}}}\delta_{-1},\qquad J\subset\{1,2,\ldots,n\}.
\end{equation}

This objective  can only be met with additional diagonalizations
of the sums.  Applying the Journ\'e Lemma as stated in Lemma~\ref{l.journe-n-1}, we can decompose $\mathcal U$ into collections
$\mathcal U_{d_1}$, for $d_1\in\mathbb N$, for which we have $2^{d_1}R\subset
V$ for $R\in\mathcal U_{d_1}$, and \begin{equation}\label{e.U_d_1BMO} \norm
\operatorname P[\mathcal U_{d_1}]b.BMO.\le{}K_{\delta_{\text{journ\'e}}}2^{(n+1)d_1}\norm
b.BMO_{-1} (\otimes_1^n\mathbb C_+).\lesssim2^{(n+1)d_1}\delta_{-1}. \end{equation} In what
follows, we shall suppress the dependence of these inequalities on
the choice of $\delta_{\text{journ\'e}}$, which comes only through
this application of Journ\'e's Lemma.   Also the (large) power of
$d_1$ is of no particular consequence.  From another part of the
estimate we can pick up a factor of $K_N2^{-Nd_1}$ for arbitrarily
large $N$.

For integers $d_2\ge{}d_1$, set \begin{gather*}
  \label{e.defXd}
\mathcal X(J,d_2):=\{(R',R)\in \mathcal W\times\mathcal U_{d_1}\mid R'\prec_J R,\ R'\subset 2^{d_2+4}R,\ R'\not\subset 2^{d_2}R\}
\\ \label{e.defZXjdd}
\mathbb X(J,d_2):=\sum_{(R',R)\in\mathcal X(J,d_2)} \overline {\ip b,v_{R'}.v_{R'}}\ip b,v_R.v_R.
\end{gather*}
In this notation, and below, we will suppress the dependence upon $d_1$, as this parameter does not directly enter
into any of the estimates.
We shall show that
\begin{equation} \label{e.ZXd<}
\norm \mathbb X(J,d_2).2.\lesssim{}2^{-d_2}\delta_{-1},\qquad J\subset\{1,2,\ldots,n\}, \ 0\le{}d_1\le{}d_2.
\end{equation}
This estimate proves (\ref{e.ZX<}).

Orthogonality enters into the estimate in the following way.
Suppose  we are given two pairs of rectangles  $(R,R') $ and $(\widetilde R,\widetilde R') $ in $\mathcal X(J)$. 
In addition, suppose  $16\abs {R'_j}<\abs{\widetilde R'_j}$ for some $j\in J$.  We conclude that the 
functions 
 $v_{R'}\overline {v_R}$ and $v_{\widetilde R'}\overline {v_{\widetilde
R}}$ are orthogonal. This is seen by examining the Fourier
supports of the wavelets. Therefore, for $|J|$-tuples of integers
$\overrightarrow\ell\in\mathbb Z^{\abs J}$, we define \begin{gather*} \label{e.defZXn}
\mathcal X(J,d_2,\overrightarrow\ell):=\{(R',R)\in\mathcal X(J,d_2)\mid
\abs{R'_j}=2^{\overrightarrow\ell_{\!j}},\ j\in{}J\}.
\\
\mathbb X(J,d_2,\overrightarrow\ell):=\sum_{(R',R)\in\mathcal X(J,d_2,\overrightarrow\ell)} \overline {\ip
b,v_{R'}.v_{R'}}\ip b,v_R.v_R. \end{gather*} Here, for simplicity, we have
assumed that $J=\{1,2,\ldots,\abs J\}$ for notational convenience.
We will continue with this assumption throughout.  All estimates
will clearly be invariant under appropriate permutation of
coordinates. In light of the orthogonality above, it is the case
that \begin{equation}\label{e.ZX->n} \norm
\mathbb X(J,d_2).2.^2\lesssim{}\sum_{\overrightarrow\ell\in\mathbb Z^{\abs J}}\norm
\mathbb X(J,d_2,\overrightarrow\ell).2.^2. \end{equation}

At this point, we can abandon orthogonality considerations
altogether in estimating this last sum.  In the sums
$\mathbb X(J,d_2,\overrightarrow\ell)$ we have the functions $v_{R'}\overline{v_R}$,
which can be dominated as \begin{align*} [\abs {R'}\abs R]^{1/2}\abs{
v_{R'}\overline{v_R}}\lesssim{}& [(\zeta_{R'}*\ind {R'})( \zeta_R*
\ind R) ]^2
\\{}\lesssim{}&\operatorname M\ind R(c(R'))^{N}\zeta_{R'}*\ind {R'},
\end{align*} 
where we take \begin{equation*} \zeta_R(x):=[1+
|x_1|\abs{R_1}^{-1}+\cdots+|x_n|\abs{R_n}^{-1}]^{-100n}. \end{equation*}
Here, $N>1$ is arbitrary, though the implied constant will
depend upon the choice of $N$.  $\operatorname M \ind R(c(R'))$ is an effective
measure of the distance between $R'$ and $R$, as $R'$ will always
have dimensions which are smaller or comparable to those of $R$.

Of course, for $(R',R)\in \mathcal X(J,d_2)$, we have $\operatorname M \ind
R(c(R'))\lesssim2^{-Nd_2}$.  Thus, 
\begin{align*}  
\norm \mathbb X(J,d_2,\overrightarrow\ell).2.
 &\lesssim{}2^{-Nd_2} \norm \widetilde
 {\mathbb X}(J,d_2,\overrightarrow\ell) .2. , \quad  \text{ where}
 \\  
 \widetilde {\mathbb X}(J,d_2,\overrightarrow\ell):={}&
 \sum_{(R',R)\in\mathcal X(J,d_2,\overrightarrow\ell)} 
 \frac{ \beta(R')}{\sqrt{\abs {R'}}} \frac{ \beta(R)}{\sqrt{\abs {R}}}
 \ind {R'}, \text{ and}\\
  \beta(R):={}&\abs{\ip b,v_R.}. 
 \end{align*} 
The top line holds
for all large integers $N$.  Thus in the argument below we can
accrue some bounded number of positive powers of $2^{d_2}$ and not
place our desired estimate in jeopardy. Hence, to obtain (\ref{e.ZXd<}) it is enough for us
to show that 
\begin{gather} \label{e.2dD} 
\sum_{\overrightarrow\ell\in\mathbb Z^{\abs J}}
\norm \widetilde{\mathbb X}(J,d_2,\overrightarrow\ell).2.^2\lesssim{}2^{8nd_2}\delta_{-1}^2,\qquad
J\subset\{1,2,\ldots,n\}.  
\end{gather}
As the terms $ \widetilde {\mathbb X}$ are sums of indicator sets of rectangles, we can appeal to 
facts about Carleson measures and, in particular, the John-Nirenberg inequalities to control these sums.

There is a final diagonalization to make.  For $|J|$-tuples of
natural numbers $\overrightarrow{d_3}\in\mathbb N^{\abs J}$, set \begin{align*}
\mathcal X(J,d_2,\overrightarrow\ell,\overrightarrow{d_3}):={}&\{(R',R)\in\mathcal X(J,d_2,\overrightarrow\ell)\mid
2^{\overrightarrow{d_3}_{j}}\abs{R'_j}=\abs{R_j},\  j\in{}J  \},
\\ 
\widetilde{\mathbb X}(J,d_2,\overrightarrow\ell,\overrightarrow{d_3}):={}&\sum_{(R',R)\in\mathcal X(J,d_2,\overrightarrow\ell,\overrightarrow{d_3})}
\frac{ { \beta(R')}}{\sqrt{\abs {R'}}} \frac{ \beta(R)}{\sqrt{\abs
{R}}} \ind {R'}
\\\nonumber
{}={}&2^{-\frac12\norm
\overrightarrow{d_3}..}\sum_{(R',R)\in\mathcal X(J,d_2,\overrightarrow\ell,\overrightarrow{d_3})} \frac{ {
\beta(R')\beta(R)}}{{\abs {R'}}} \ind {R'}, \end{align*} where  $\norm
{\overrightarrow{d_3}}..:=\sum_{j=1}^{\abs J} \abs{\overrightarrow{d_3}_j}$. The leading
term of the last line suggests that indeed $\overrightarrow{d_3}$ is a
diagonalization parameter.

With the relative sizes of $R'$ and $R$ fixed by the choice of $J\subset\{1,2,\ldots,n\}$ and by the choice of $\overrightarrow{d_3}$,
observe that for each $R'$, there can be at most $O(2^{n d_2})$ possible choices of $R$ so that 
$(R,R')\in \mathcal X(J,d_2,\overrightarrow\ell,\overrightarrow{d_3})$.     Again, we can afford to lose some bounded number of powers of $2^{d_2} $ in 
our estimates. We take
\begin{align*}\label{e.devzoY}
\mathcal Y(J,d_2,\overrightarrow\ell,\overrightarrow{d_3}):={}&\{R' \mid (R',R)\in\mathcal X(J,d_2,\overrightarrow\ell,\overrightarrow{d_3})\text{${}$ for some $R\in\mathcal R$}\},
\\
\mathcal Y(J,d_2,\overrightarrow{d_3}):={}&\bigcup_{\overrightarrow\ell\in\mathbb Z^{\abs
J}}\mathcal Y(J,d_2,\overrightarrow\ell,\overrightarrow{d_3}). \end{align*} 
Let $\pi \,:\,
\mathcal Y(J,d_2,\overrightarrow\ell,\overrightarrow{d_3})\longrightarrow\mathcal R$ be such that
$(R',\pi(R'))\in\mathcal X(J,d_2,\overrightarrow\ell,\overrightarrow{d_3})$. We would need to consider  $O(2^{n d_2})$ possible choices for 
this function.  Below, we will consider just some arbitrary choice of this function $\pi$, and then merely sum over the 
possible choices of $\pi$, accruing a harmless term of  $O(2^{n d_2})$. Set 
\begin{equation} 
\label{e.defzy}
		\begin{split} 
\mathbb Y(J,d_2,\overrightarrow\ell,\overrightarrow{d_3}):={}&\sum_{R'\in\mathcal Y(J,d_2,\overrightarrow\ell,\overrightarrow{d_3})}
\frac{  \beta(R')}{\sqrt{\abs {R'}}} \frac{ \beta(\pi(R'))}{\sqrt{\abs
{\pi(R')}}} \ind {R'}
\\{}={}& 2^{-\frac12\norm \overrightarrow{d_3}..}
\sum_{R'\in\mathcal Y(J,d_2,\overrightarrow\ell,\overrightarrow{d_3})}  \frac{
\beta(R')\beta(\pi(R'))}{{\abs {R'}}} \ind {R'}. 
	\end{split}
	\end{equation}

The specific estimate we prove  is:
\begin{equation}\label{e.2do} 
\begin{split}
\norm
\mathbb Y(J,d_2,\overrightarrow{d_3}).2.^2:={}&\sum_{\overrightarrow\ell\in\mathbb Z^{\abs J}} \norm
\mathbb Y(J,d_2,\overrightarrow\ell,\overrightarrow{d_3}).2.^2
\\{}\lesssim{}&2^{8nd_2-\frac14\norm \overrightarrow{d_3}..}\delta_{-1}^2,\qquad J\subset\{1,2,\ldots,n\}.  
\end{split} 
\end{equation} This is summed over $\overrightarrow{d_3}\in\mathbb N^{\abs J}$ to prove (\ref{e.2dD}),  and so will complete our proof. The proof of this
inequality is taken up in the next subsection. We achieve an
exponential decay in parameters $d_1,\ d_2,\ \overrightarrow{d_3}$.

We shall rely repeatedly  on the estimates \begin{align}\label{e.use1}
\sum_{R\in \mathcal U_{d_1}}\beta(R)^2&\lesssim{}2^{2(n+1)d_1}\delta_{-1}^2,
\\ \label{e.use2}
\sum_{R'\in\mathcal Y(J,d_2,\overrightarrow{d_3})}\beta(R')^2&\lesssim{} 2^{2nd_2}, \qquad J\subset\{1,2,\ldots,n\},\
 \overrightarrow{d_3}\in\mathbb N^{\abs J}.
\end{align}
The first of these has the critical gain by a factor of $\delta_{-1}^2$, as 
follows from (\ref{e.U_d_1BMO}).  The
second estimate follows from the fact that $b$ is in $BMO$ and that the rectangles $R'$ in $\mathcal Y(J,d_2,\overrightarrow{d_3})$
 are contained in $\{ \operatorname M\ind {\sh {\mathcal U} }>c2^{-nd_2}\}$ since
 $R'\subset 2^{d_2+4}R$ for some $R\in \mathcal  U_{d_1}$.  Here $\operatorname M$ is the strong maximal function.

At this point we recap the notations.
\begin{itemize}
\item{$\delta_{-1}=\norm b.BMO_{-1} (\otimes_1^n\mathbb C_+).$.}

\item{$d_1$ is associated to the measure of embeddedness of rectangles $R\in \mathcal U$.}

\item{$d_2$ is a (crude) measure of the separation between the rectangles $R'\in\mathcal W$ and $R\in\mathcal U_{d_1}$ for  $d_2\ge{}d_1$.}

\item{$J\subset\{1,2,\ldots,n\}$ is that set of coordinates for which one has some orthogonality.}

\item{$\overrightarrow\ell\in\mathbb Z^{\abs J}$ specifies the side lengths of $R'$ for those coordinates  $j\in J$.}

\item{$\overrightarrow{d_3}\in\mathbb N^{\abs J}$ specifies how much bigger $R$ is than $R'$ in the coordinates $j\in J$. }

\item{$(R',\pi(R'))\in\mathcal X(J,d_2,\overrightarrow\ell,\overrightarrow{d_3})$ and $\abs{\pi(R')_j}=2^{{\overrightarrow{d_3}}_j}\abs{R'_j}$ for $j\in J $, otherwise for $j\not \in J$,
$\abs{\pi(R')_j} {}\simeq{}\abs{R'_j}$. }

\item $   \beta(R):=\abs{\ip b,v_R.} $.
\end{itemize}

\subsection*{The Bounds for the Paraproducts}

The argument varies depending upon the cardinality of
$J\subset\{1,2,\ldots,n\}$.   While we can formalize issues in a
way that is uniform with respect to $\abs J$, we present four
subsections, to emphasize the differences that come about due to
the increasing number of parameters.

\subsubsection*{The Case of $J=\{1,2,\ldots,n\}$.}

The key point is that the sum in (\ref{e.2do}) simplifies
considerably, as all the side lengths of $R'$ are specified by the
parameter $\overrightarrow\ell$.    In particular, the rectangles $R' $ occurring in the sum 
in (\ref{e.defzy}) are pairwise disjoint.  
Thus, the $L^2$ norm in (\ref{e.2do}) will
simplify to 
\begin{align*}
\norm  \mathbb Y(J,d_2,\overrightarrow\ell,\overrightarrow{d_3}).2.^2
&{}={}\sum_{R'\in\mathcal Y(J,d_2,\overrightarrow\ell,\overrightarrow{d_3})}
\frac{\beta(R')^2 \beta(\pi(R'))^2}{\pi(|R'|)} 
\\&{}={}2^{-\norm \overrightarrow{d_3}..}\sum_{R'\in\mathcal Y(J,d_2,\overrightarrow\ell,\overrightarrow{d_3})}
\frac{\beta(R')^2 \beta(\pi(R'))^2}{|R'|} 
\\&{}\le{}2^{-\norm \overrightarrow{d_3}..} \sup_{R' }\frac{ \beta(R')^2}{\abs{R' } }
\sum_{R'\in\mathcal Y(J,d_2,\overrightarrow\ell,\overrightarrow{d_3})}\beta(\pi(R'))^2.
\end{align*} 
As $b $ has $BMO $ norm one, the supremum above is bounded by $1 $.  Then sum over $\overrightarrow\ell\in\mathbb Z^{\abs J}$ and use  (\ref{e.use1}) to see that 
\begin{equation*} \norm
\mathbb Y(J,d_2,\overrightarrow{d_3}).2.^2\lesssim{}\delta_{-1}^2 2^{3n
d_2-\norm \overrightarrow{d_3}..}. 
\end{equation*}
Recall that we can tolerate a few positive powers of $ d_2 $.  This case is complete.


\subsubsection*{ The Case of $J=\{1,2,\ldots,n - 1\}$.}
Now, the rectangles $R'\in\mathcal Y(J,d_2,\overrightarrow\ell,\overrightarrow{d_3})$ are only
permitted to vary in the last coordinate. That is, the
corresponding sums are as complex as those of one parameter
Carleson measures. So we can explicitly compute \begin{align*} \norm
\mathbb Y(J,d_2,\overrightarrow\ell,\overrightarrow{d_3}).2.^2={} &2^{-\norm \overrightarrow{d_3}..}
\sum_{R'\in\mathbb Y(J,d_2,\overrightarrow\ell,\overrightarrow{d_3})}
    \frac{  \beta(R')\beta(\pi(R'))}{{\abs {R'}}}
    \\&\qquad {}\times{}
    \sum_{\substack{R''\in\mathbb Y(J,d_2,\overrightarrow\ell,\overrightarrow{d_3})\\ R''\subset R'}}
    {  \beta(R'')} { \beta(\pi(R''))}.
\end{align*}
With the specific way the innermost sum is formed, observe that
\begin{align*}
\sum_{\substack{R''\in\mathbb Y(J,d_2,\overrightarrow\ell,\overrightarrow{d_3})\\ R''\subset R'}}\beta(R'')^2&\lesssim{}\delta_{-1}^2\abs {R'},
\\
\sum_{\substack{R''\in\mathbb Y(J,d_2,\overrightarrow\ell,\overrightarrow{d_3})\\ R''\subset
R'}}\beta(\pi(R''))^2&\lesssim{} 2^{n d_2+\norm
\overrightarrow{d_3}..}\delta_{-1}^2\abs {R'}. \end{align*} The first estimate is
obvious, while the second estimate follows from the fact that the
rectangles $\pi(R')$ are contained in \begin{equation*} \otimes_{j=1}^n
2^{d_2+\overrightarrow{d_3}_j}R'_j . \end{equation*} In this last display, set the last
coordinate of $\overrightarrow{d_3}$ to be zero. Applying these observations,
Cauchy-Schwarz,  (\ref{e.use2}), and $\norm b.BMO (\otimes_1^n\mathbb C_+).=1$ we see
that \begin{align*} \norm \mathbb Y(J,d_2,\overrightarrow{d_3}).2.^2\lesssim{}&\delta_{-1}^2
2^{-\frac12\norm \overrightarrow{d_3}..}2^{\frac{nd_2}{2}}
\sum_{R'\in\mathbb Y(J,d_2,\overrightarrow\ell,\overrightarrow{d_3})} \beta(R')\beta(\pi(R'))
\\{}\lesssim{}&
\delta_{-1}^22^{\frac52 nd_2-\frac12\norm \overrightarrow{d_3}..}. \end{align*} This
completes this case.


\subsubsection*{The Case of $0<\abs J<n - 1$.}
The argument in this case could be adapted to treat the general case.
We would like to indicate the additional difficulty that one faces in this case.
 The side lengths of $R'$ are fixed for those coordinates in $J$, and  completely specified by $\overrightarrow\ell\in\mathbb Z^{\abs J}$.
 The remaining
 side lengths of $R'$ are then permitted to vary.  Thus, the ways that two possible choices of
 $R',R''\in\mathcal Y(J,d_2,\overrightarrow\ell,\overrightarrow{d_3})$ can intersect are as general as the intersections
 of two dyadic rectangles of dimension $n-\abs J$.

 Nevertheless, one can implement a method of proof that follows the lines of the case $J=\{1,2,\ldots,n - 1\}$,
 provided one takes advantage of the John-Nirenberg inequality, which we now state in the form used. 
 For rectangles $R \in \mathcal R$ and non-negative constants $a_R$ 
 for which $\sum_{R \subset W} a_R \le |W|$  for all open sets $W \subset \mathbb R^n$, we have
 $$
 \NORm \sum_{R\subset W}\frac{a_R}{|R|} \ind{R} .p.\lesssim|W|^{1/p}, \qquad 1<p<\infty.
 $$
 We use this to obtain  the following extensions of the
 inequalities (\ref{e.use1}) and (\ref{e.use2}).  In the first place, we have
 \begin{equation} \label{e.Use1}
 \NORm \Biggl[ \sum_{R'\in \mathcal Y(J,d_2,\overrightarrow{d_3})} \frac{\beta(R')^2 } {{\abs {R'}}} \ind {R'}
 \Biggr]^{1/2} .p.\lesssim{} 2^{2nd_2},\qquad 1<p<\infty.
 \end{equation}
 This is available to us from the fact that $b $ is in $BMO $ with norm one.
 A similar fact is
 \begin{equation} \label{e.Use2}
 \NORm \Biggl[ \sum_{R'\in \mathcal Y(J,d_2,\overrightarrow{d_3})} \frac{\beta(\pi(R'))^2 } {{\abs {\pi(R')}}} \ind {R'}
 \Biggr]^{1/2} .p.\lesssim{}\delta_{-1} 2^{2nd_2},\qquad 1<p<\infty.
 \end{equation}
The important  features of these estimates are that they are  independent of
$\overrightarrow{d_3}$,  uniform in $\overrightarrow\ell\in\mathbb Z^{\abs J}$, and in the second estimate we have the gain of $\delta_{-1}$. 

Now, the second estimate does not immediately follow from a $BMO $ estimate, due to the fact that we have a 
mismatch between $R' $ and $\pi(R') $ in (\ref{e.Use2}). 
Due to the
John-Nirenberg inequality, (\ref{e.Use2})  will follow from the
estimate \begin{equation*} \sum_{\substack{R'\subset W\\
R'\in\mathcal Y(J,d_2,\overrightarrow{d_3})}} \beta(\pi(R'))^2\lesssim{}
2^{4nd_2+\norm \overrightarrow{d_3}..}\delta_{-1}^2 \abs {W},\qquad W\subset
\mathbb R^n. \end{equation*} As this estimate is uniform in the choice of $W$, it
provides a bound for a  Carleson measure to which the
John-Nirenberg inequality applies.
 For a given $W\subset\mathbb R^n$, it is the case that for all rectangles $R'$ that contribute to this sum,
 the rectangle $\pi(R')$ is contained in a set which is given in the first place by a strong maximal
 function applied to $W$.  Set
\begin{equation*} W_0:=\{\operatorname  M\ind W \ge{}c2^{-n d_2} \}, \end{equation*} for an appropriate
choice of $c$.   This set, so constructed, will contain a
translation of $R'$ which is contained in $\pi(R')$.  The point to
keep in mind is that $\pi(R')$ is $2^{\overrightarrow{d_3}_j}$ times longer than $R'$ in the
coordinate $j\in J$.  Thus, in that coordinate, we should
apply a one dimensional maximal function with threshold
$2^{-\overrightarrow{d_3}_j}$.  Namely, for $j\in{}J=\{1,2,\ldots,\abs J\}$,
we inductively define \begin{equation*} W_j:=\{ \operatorname M_j\ind
{W_{j-1}}>c2^{-\overrightarrow{d_3}_j}\}. \end{equation*} For appropriate constant $c$,
we will have $\pi(R')\subset W_{{\abs J}}$.  And we certainly
have $\abs {W_{\abs J}}\lesssim{} 2^{\norm \overrightarrow{d_3}..+2nd_2}\abs
W$.  This completes the proof of (\ref{e.Use2}).

Estimates  (\ref{e.Use1}) and (\ref{e.Use2}) are not in themselves enough to complete the
proof, as there is no decay in the quantity $\norm \overrightarrow{d_3}..$.
But, they do show that 
\begin{align} \nonumber&\NOrm \Bigl[\sum_{\overrightarrow\ell\in\mathbb Z^{\abs
J}} \mathbb Y(J,d_2,\overrightarrow\ell,\overrightarrow{d_3})^2\Bigr]^{1/2}
.4.\\ \label{e.Use3}
\lesssim{}&\NORm \Biggl[ \sum_{R'\in \mathcal Y(J,d_2,\overrightarrow{d_3})}
\frac{\beta(R')^2 } {{\abs {R'}}} \ind {R'}
 \Biggr]^{1/2} .8.\NORm \Biggl[ \sum_{R'\in \mathcal Y(J,d_2,\overrightarrow{d_3})} \frac{\beta(\pi(R'))^2 } {{\abs {\pi(R')}}} \ind {R'}
 \Biggr]^{1/2} .8.\\\nonumber
 \lesssim{}&2^{4nd_2}\delta_{-1} . \end{align} Namely, we have an estimate on
the $L^4$ norm that is uniform with respect to $\overrightarrow{d_3}$. This
will permit us to select a set which decays with respect to this
parameter. On this set, we will not attempt to estimate the $L^2$
norm in (\ref{e.2do}).
The set we take is
\begin{equation*}
E:=\bigcup_{\overrightarrow\ell\in\mathbb Z^{\abs J} }\{\operatorname  M 
\mathbb Y(J,d_2,\overrightarrow\ell,\overrightarrow{d_3})
>\delta_{-1}2^{\frac18\norm \overrightarrow{d_3}..} \}.
\end{equation*}
Here, we use the strong maximal function $ \operatorname M$.
This set has measure $\abs E\lesssim{}2^{16nd_2-\frac12\norm \overrightarrow{d_3}..}$, due to the large $L^p$ norms we have in  (\ref{e.Use3}).

To complete the argument in this case, it suffices to show that
\begin{equation*}
\sum_{\overrightarrow\ell\in\mathbb Z^{\abs J} } \int_{\mathbb R^n-E} \abs{ \mathbb Y(J,d_2,\overrightarrow\ell,\overrightarrow{d_3})}^2 \; dx
\lesssim\delta_{-1}2^{3nd_2-\frac14\norm \overrightarrow{d_3} ..}.
\end{equation*}
We will expand the square on the left hand side.   Integrating $ \mathbb Y(J,d_2,\overrightarrow\ell,\overrightarrow{d_3}) $ over $R'-E $, we will 
lose  a factor of $2^{\frac18\norm \overrightarrow{d_3}..}$.  But from $\abs{ \pi(R')}=2^{\norm \overrightarrow{d_3}..}\abs{R'} $ we will gain a factor 
of $2^{\frac12\norm \overrightarrow{d_3}..}$.  Specifically,
\begin{align*}
 \sum_{\overrightarrow\ell\in\mathbb Z^{\abs J}}\int_{\mathbb R^n-E} \abs{ \mathbb Y(J,d_2,\overrightarrow\ell,\overrightarrow{d_3})}^2 \; dx\le{}&
 \sum_{\overrightarrow\ell\in\mathbb Z^{\abs J}}\sum_{R'\in \mathcal Y(J,d_2,\overrightarrow\ell,\overrightarrow{d_3})}
\frac{\beta(R')}{\sqrt{\abs {R'}}}\frac {\beta(\pi(R'))}{\sqrt{\abs{\pi(R')}}}\times
\\{}&\qquad
\int_{R'-E}\mathbb Y(J,d_2,\overrightarrow\ell,\overrightarrow{d_3}) \; dx
\\{}\le{}&\delta_{-1} 2^{(\frac18-\frac12){\norm \overrightarrow{d_3}..}} \sum_{\overrightarrow\ell\in\mathbb Z^{\abs J}}\sum_{R'\in \mathcal Y(J,d_2,\overrightarrow\ell,\overrightarrow{d_3})} \beta(R')\beta(\pi(R'))
\\{}\lesssim{}&\delta_{-1}^22^{(2n+1)d_2-\frac14\norm \overrightarrow{d_3} ..}.
\end{align*}
This estimate follows from the definition of the set $E$ and (\ref{e.use1}) and (\ref{e.use2}).


\subsubsection*{The Case of $J=\emptyset$.}
In this case, both $\overrightarrow\ell$ and $\overrightarrow{d_3}$ are not present, and
the rectangles $R'$ and $\pi(R')$ have comparable lengths in all
coordinates.    But we do have (\ref{e.Use1}) and (\ref{e.Use2}),
and they directly prove the desired estimate
 \begin{align*} \norm
\mathbb Y(\emptyset,d_2).2.\lesssim{} & \NORm \Biggl[ \sum_{R'\in
\mathcal Y(\emptyset,d_2)} \frac{\beta(R')^2 } {{\abs {R'}}} \ind {R'}
 \Biggr]^{1/2} .4.
 \NORm \Biggl[ \sum_{R'\in \mathcal Y(\emptyset,d_2)} \frac{\beta(\pi(R'))^2 } {{\abs {\pi(R')}}} \ind {R'}
 \Biggr]^{1/2} .4.
\\{}\lesssim{}&\delta_{-1} 2^{4nd_2}.
\end{align*}
And this completes this case.



\section{Journ\'e's Lemma}\label{s.journe}

We state a version of the Lemma of Journ\'e \cite{MR87g:42028} that is implicit in Pipher's extension \cite{MR88a:42019}, and  interfaces well with our notion of
a restricted $BMO$ norm, namely $BMO_{-1} (\otimes_1^n\mathbb C_+)$.   We first state the lemma in a purely
geometric fashion, and then return to a formulation that is more specific to our needs in this paper.

\subsection*{The Geometric Formulation}

Given a   collection $\mathcal U$ of   dyadic
rectangles whose shadow has finite area,    suppose that $V\supset \sh{\mathcal U}$. For  rectangles $R\in\mathcal U$,  define
\begin{equation*}
\emb R.V.:=\sup\{\mu\ge1\mid \mu \Rho 1\times\Rho 2\times \cdots\times \Rho n\subset V\}.
\end{equation*}
For an arbitrary subset $\mathcal U'\subset \mathcal U$, let
\begin{equation*}
F(I,k, \mathcal U'):=\bigcup\{ I\times R'\mid I\times R'\in \mathcal U',\ 2^{k-1}\le\emb  I\times R'.V.<2^k\}.
\end{equation*}

\begin{lemma}\label{l.few-small}
 For all $\delta,\epsilon>0$,  we can select $V\supset \sh {\mathcal U} $
with $\abs V\le{}(1+\delta)\abs {\sh {\mathcal U} }$, for which   we have the uniform estimate\footnote{We have stated the lemma in the formulation for the first
coordinate to ease the burden of notation. In application, we will use this in an arbitrary choice of coordinate.}
\begin{equation}  \label{e.few-small}
\sum_{k=1}^\infty \sum_{I\in\mathcal D} 2^{-\epsilon k}\abs{ F(I,k,\mathcal U')}\lesssim{}\abs{\sh {\mathcal U'}},\qquad \mathcal U'\subset\mathcal U.
\end{equation}
The implied constants in these estimates depend only on dimension and the choices of $\epsilon,\delta$.
\end{lemma}

  Consider a collection of rectangles $\mathcal U$ in which all the first coordinates are the same.  Then the embeddedness
  is necessarily of order $1$.  This example shows that the lemma above
  must be formulated in this fashion.

We begin the proof with a careful description of how to select the set $V$.  If we were not too concerned about the upper bound
on the measure of $V$, in other words if the bound 
$\abs V\lesssim{}\abs {\sh {\mathcal U}}$ were enough, then we could simply
take $V=\{\operatorname M\ind {\sh {\mathcal U} }>\frac12\}$.   
For our needs, however, this choice of $V$ is completely inappropriate.

We  need the notion of {\em shifted dyadic grids}, which is a modification of
an observation due to M.~Christ  defined as follows.
The definition of the grids depends upon a choice of integer $\mathsf d$, and we will set  $\delta=(2^{\mathsf d}+1)^{-1}$.   For integers
$0\le{}b<{}\mathsf d$, and $\alpha\in\{\pm(2^{\mathsf d}+1)^{-1}\}$, let
 \begin{gather} \label{e.shifted-grids}
 \mathcal D_{\mathsf d,b,\alpha}:=\{2^{k\mathsf d+b}((0,1)+j+(-1)^k\alpha)\mid k\in\mathbb Z,\ j\in\mathbb Z\},
 \\\nonumber
 \mathcal D_{\mathsf d}:=\bigcup_\alpha\bigcup_{b=0}^{\mathsf d-1}\mathcal D_{\mathsf d,b,\alpha}.
 \end{gather}
 One checks that  $\mathcal D_{\mathsf d,b,\alpha}$ is a grid.  Indeed, it suffices to assume
 $\alpha=(2^{\mathsf  d}+1)^{-1}$ and that $b=0$.  Checking the grid structure can be done by
 induction.  And it suffices to check that the intervals in $\mathcal D_{\mathsf d,0,\alpha}$ of
 length one are a union of intervals in $\mathcal D_{\mathsf d,0,\alpha}$ of length $2^{-\mathsf d}$.
 One need only check this for the interval $(0,1)+\alpha$.  But certainly
  \begin{align*}
 (0,1)+\frac1{(2^{\mathsf d}+1)}
 {}&{}=\bigcup_{j=0}^{2^{\mathsf d}-1} (0,2^{-\mathsf d})+\frac j{2^{\mathsf d}}+\frac1{(2^{\mathsf d}+1)}
 \\&{}={}
 \bigcup_{j=0}^{2^{\mathsf d}-1} (0,2^{-\mathsf d})+\frac {j+1}{2^{\mathsf d}}-\frac {1}{2^{\mathsf d}(2^{\mathsf d}+1)}.
 \end{align*}

   What is more important concerns the collections
 $\mathcal D_{\mathsf d}$.  For each dyadic interval
 $I\in\mathcal D$, $I\pm\delta\abs I\in\mathcal D_{\mathsf d}$.\footnote{The problem we are avoiding here is that the dyadic grid distinguishes  dyadic rational points. At the point $0$ for instance,
 $((0,1)-\delta) \not\subset(0,2^k)$ for all integers $k$,
 regardless of how big $k$ is.}   Moreover, the maximal function $\operatorname M^{\mathcal D_{\mathsf d}}$ maps $L^1$ into
 $L^{1,\infty}$\footnote{In fact, taking $\mathsf d=1$, it is routine to check that $M^{\mathcal D_{\mathsf 1}}$ dominates
an absolute multiple of the usual maximal function,  thus, proving that it satisfies the weak type inequality.}
 with norm at most $2\mathsf d\simeq\abs{\log \delta}$. In fact we need the finer estimate, valid for all 
 choices of $0<\delta<1 $ and integers $\mathsf d $,
\begin{equation*}
\abs{ \{\operatorname M^{\mathcal D_{\mathsf d}}\ind
U>1-\delta\}}\le{}(1+K\delta\mathsf d)\abs U 
\end{equation*} 
for all subsets $U$ of
the real line of finite measure and some constant $K$. 
This will be an 
effective estimate since the   value of $\mathsf d $ we will consider is 
$\mathsf d\simeq\lceil -\log_2 \delta\rceil $.   
 To see this estimate,  note that 
\begin{align*} 
\abs{\{ \operatorname M^{\mathcal D_{\mathsf d}}\ind U
>1-\delta\}}\le{}& \abs{U}+\sum_{b=0}^{\mathsf d-1}\sum_{\alpha\in \{\pm
(2^{\mathsf d}+1)^{-1}\}} \abs{ U^c\cap \{ M^{\mathcal D_{\mathsf d, b,
\alpha}}\ind U >1-\delta\}}
\\
{}\le{}&(1+2\mathsf d[(1-\delta)^{-1}-1])\abs U.
\end{align*}

The main line of the argument can now begin.
We take
$\delta=(1+2^{\mathsf d})^{-1}$ for an integer $\mathsf d$. We use the maximal functions
$\operatorname M^{\mathcal D_{\mathsf d}}$, but only in the last step of the induction. 
Initialize  $\text{Enl}(n+1,\mathcal U'):=\sh {\mathcal U' }$, so that we use backwards induction.   Inductively define
\begin{equation}  \label{e.V-def} 
	\begin{split} 
\text{Enl}  (j,\mathcal U'){}&:={}
\{\operatorname  M^{\mathcal D}_j\ind {\text{Enl}(j+1,\mathcal U') }>1-\delta ^{2^j}\},\qquad n\ge{} j\ge2,
\\
V&{}:=\{\operatorname  M^{\mathcal D_{\mathsf d}}_1\ind {\text{Enl} (2,\mathcal U')}>1-\frac\delta2\},
	\end{split}
	\end{equation}
where the subscript on the maximal functions denotes the  coordinate in which the maximal
function 
is applied.  Then it is the case that 
$\abs V\le{}(1+K\delta\abs{\log \delta})\abs {\sh{\mathcal U'} }$,
where the constant $K$ depends only on the dimension $n$.

 Now pass to a further subset  $\mathcal U''\subset \mathcal U'$  such that  for all $R, R'\in\mathcal U''$, we have  
 $2^{k-1}\le\emb R.V.\le2^k $   and   if it is the
case that $\abs {R'_1}<\abs {R_1}$, then we have the stronger inequality $40\cdot 2^k \delta^{-1}\abs
{R'_1}<\abs {R_1}$.  We term this assumption ``separation of
scales''   in the first coordinate.  Under these assumptions, estimate (\ref{e.few-small})
reduces to \begin{equation}  \label{e.frozen} 
\sum_{I\in\mathcal D} \abs{
F(I,k,\mathcal U'')}\lesssim{}\abs{\sh {\mathcal U'}},\qquad \mathcal U'\subset\mathcal U.
\end{equation} 
Sufficiency is seen by noting that  obtaining separation of
scales necessitates  dividing the rectangles into approximately
$\log 2^k\delta^{-1}$ subclasses.  Then multiplying by
$2^{-\epsilon k}$, one is able to sum over all scales $k$.

Our strategy is to define, for dyadic intervals $I$, sets $H(I)$ that are disjoint in $I$, contained in
$\sh {\mathcal U'}$, and for which
\begin{equation*}
\ind {F(I ,k,\mathcal U'')}\lesssim{} \operatorname M \ind {H(I)}
\end{equation*}
for an appropriate maximal function $\operatorname M$, where the implied constant is permitted to depend upon 
$\delta $ and dimension $n $.  It will in fact be of the order $O(\delta ^{-2^{n+1} }) $.
An appeal to the Fefferman-Stein maximal
inequalities \cite{MR44:2026} will then prove (\ref{e.frozen}).

The sets $H(I) $ are defined to be $F(I ,k,\mathcal U'')-G(I) $, where 
\begin{equation*}
G(I):=\bigcup_{\substack{ I'\in\mathcal
D\\I\subset_{\not=}I'}}F( I',k,\mathcal U''). 
\end{equation*} 
By our separation of scales,  
the minimal dyadic interval  $\widetilde I $   that contributes to this union  contains $I $ and  satisfies  
\begin{equation*}{}
40\cdot 2^k \delta^{-1}\abs I<\abs {\widetilde I }\le{}80\cdot 2^k \delta^{-1}\abs I .
\end{equation*}

We need to show that if $R  $ is a dyadic rectangle with 
$\abs {R\cap G(I) }>(1-\delta ^{2 ^{n+1 } })\abs R $, then $\emb R.V.\ge{}2^{k+1} $, and hence it can't be  among those 
rectangles that contribute to $ F(I,k,\mathcal U'') $.\footnote{If we could use the  strong maximal function to define embeddedness, 
this conclusion would be immediate. Our more subtle definition of embeddedness appears to force the more complicated 
argument that follows.}
This will be accomplished by the following device.  We will show that 
\begin{equation} \label{e.will-do} 
\abs{ \widetilde I \times R_{2}\times\cdots\times R_{n }\cap\text{Enl}  (2,\mathcal U') }\ge{} 
(1-\delta^{2^2})\abs{ \widetilde I \times R_{2}\times\cdots\times R_{n} }.
\end{equation}
As $\widetilde I\pm\frac\delta2\abs {\widetilde I }\in \mathcal D_d $, and this is the grid we use in the final 
step in the construction of $V$, we conclude that the rectangle  $\widetilde I \times R_{2}\times\cdots\times R_{n }$ is inside of $V $.  
Even $\delta\abs{\widetilde I}>20\cdot2^k 
\abs I $, therefore we see that $\emb R.V.\ge{}2^{k+1} $, as desired.

We turn to the proof of (\ref{e.will-do}).  The sequence of powers of $\delta$ that appear in the definition of $V $, (\ref{e.V-def}), 
is explained in part by the next proposition. 

\begin{proposition}\label{p.probability}
  Let $0\le{} X\le1 $ be a random variable on a probability space satisfying $\mathbb E X=1-\eta $. Then, 
\begin{equation*}
\mathbb P(X<1-\sqrt \eta)\le\sqrt \eta.
\end{equation*}
\end{proposition} 

\begin{proof} 
Setting  $p=P(X<1-\sqrt \eta)$,  the inequality 
\begin{equation*}
1-\eta\le{}\mathbb E X\le(1-\sqrt \eta)p+1-p\,
\end{equation*}
will prove the proposition.
\end{proof}

We continue with the  language of probability.  Let $(\Omega_j,\mathbb P_j) $ be the probability spaces 
\begin{align*}
\Omega_1&:=\widetilde I,\\
\Omega_j&:=\widetilde I\times R_{2}\times\cdots\times R_{j},\qquad 2\le{}j\le{}n,
\end{align*}
and let $\mathbb P_j$ be normalized Lebesgue measure on $\Omega_j $.  The first of the relevant sequence of random variables 
on these spaces is 
\begin{align*}
X_{n-1}(x)&{}:=\frac {\abs{(\{x\}\times R_{n})\cap G(I)} } {\abs{R_{n} } },\qquad x\in \Omega_{n-1},
\\
A_{n-1}&{}:=\{ x\mid X_{n-1}(x)>1-\delta^{2^n}\} .
\end{align*}
Since $\abs {R\cap G(I) }>(1-\delta ^{2 ^{n+1 } })\abs R $,  $\mathbb E X_{n-1}\ge{}1-\delta^{2^{n+1 }} $, and
applying the proposition, $\mathbb P _{n-1} (A_{n-1})>{}
1-\delta^{2^n}$.  Continuing by reverse induction,  define for $n-1\ge j\ge 2$
\begin{align*}
X_{j-1}(x)&{}:=\frac { \abs{( \{x\}\times R_{j})\cap A_j } } {\abs{ R_{j} }} ,\qquad x\in \Omega_{j-1},\\
A_{j-1}&{}:=\{x\mid X_{j-1}(x)>1-\delta^{2^{j} }\} .
\end{align*}
Induction gives us the conclusion that $ \mathbb P_1(A_1)>1-\delta^{2^2} $.  This implies 
(\ref{e.will-do}) by inspection of definitions and so completes the proof.


\subsection*{A Second Geometric Formulation}

We need a certain variant of the previous lemma.
Given a set of rectangles $\mathcal U$, we let $\Emb\cdot..\,:\,\mathcal U\to[1,\infty)$ be a map
from $\mathcal U$ to the reals greater than one.  And we take $\imath\,:\, \mathcal U\to\{1,2,\ldots,n\}$
which is simply a choice of coordinates.  Based on these two data, for any subset $\mathcal U'\subset\mathcal U$ we set
\begin{equation}  \label{e.F}
F(I,k,m,\mathcal U')={}\bigcup\{ R\in\mathcal U'\mid 2^k\le{}\Emb R..<{}2^{k+1},\ \imath(R)=m,\ \Rho m=I\}.
\end{equation}

\begin{lemma}\label{l.journe}
  Fix $\delta,\,\epsilon>0$.  For any collection of rectangles $\mathcal U$ with finite shadow, we can select
$V\supset\sh {\mathcal U}$, and data $\Emb \cdot..$ and $\imath$ so that the following conditions hold.
\begin{gather*}
\abs V\le{}(1+\delta)\abs{\sh {\mathcal U}},
\\
\Emb R..R\subset V,\qquad R\in\mathcal U,
\\
\sum_{k=0}^\infty\sum_{m=1}^n\sum_{I\in \mathcal D}2^{-(n+\epsilon) k}\abs{ F(I,k,m,\mathcal U')}\lesssim{}\abs{\sh {\mathcal U'}}.
\end{gather*}
The implied constant in the last line depends only on $\delta,\,\epsilon$ and dimension.
\end{lemma}

The essential points for us are that the set $V$ is not much larger than the shadow of $\mathcal U$, and that the rectangles
$R\in\mathcal U$, after a dilation {\em uniform in all coordinates} by the embeddedness quantity,  is contained in the
enlarged set $V$.    We find that the embeddedness quantity in the last line requires a large negative power, but that
is a completely harmless fact in the context of the application we have in mind.

Again, the fact that the dyadic intervals distinguish certain points causes some difficulties for
us, and we appeal to the
shifted dyadic intervals  (\ref{e.shifted-grids}) of the previous subsection,
though our needs are not so refined in the current context.  Let $\mathcal S$ be the union of the
dyadic intervals with the two collections  $\mathcal D_{\pm}:=\mathcal D_{\mathsf 1,0,\pm\frac13}$.  For any interval $I$ of the real line,
we can find an interval $J\in\mathcal S$ with $I\subset J\subset4I$.  Indeed, let $J'$ be the maximal
dyadic interval contained in $4I$ with $\abs{I\cap J'}\ge\frac12\abs I$. If $I\subset J'$ we
are done, so assume that this is not the case.  We necessarily have $\abs{J'}\ge\frac94\abs I$,
so that one of the two intervals $J'\pm\frac13\abs{J'}$ contains $I$.  Both of these intervals are in $\mathcal S$, so we are done.

The method of proof requires that we  apply Lemma~\ref{l.few-small}, although
we find it necessary to apply it both inductively and to a wide range of possible collections of
rectangles.  In fact, it is useful to us that this lemma applies not just to collections of
a subset  $\mathcal U\subset \otimes_{j=1}^n\mathcal S$ such that the shadow of $\mathcal U$ is of finite measure.   
It also applies to all possible subsets of $\mathcal U$.

We apply Lemma~\ref{l.few-small}  to $\mathcal U^0:=\mathcal U$.  Thus, we get a set $V^{1}\supset\sh {\mathcal U^0}$, with $\abs{ V^1}\le{}(1+\delta)\abs {\sh {\mathcal U^0}}$, so that
for
\begin{equation*}
\operatorname{emb}^1 (R,V^1):=\sup\{\mu\ge1\mid \mu \Rho 1\times\Rho 2\times \cdots\times \Rho n\subset V^1\}.
\end{equation*}
we have the conclusion of Lemma~\ref{l.few-small} holding.
We then construct $\mathcal U^1\subset\otimes_{j=1}^n\mathcal S$. Set
\begin{align*}
\mathcal U^1:=\{ \gamma\times \otimes_{j=2}^n \Rho j\mid& R\in\mathcal U,\
    \gamma\in\mathcal S,\\ {}&\qquad (\Rho 1\cup \tfrac14\operatorname{emb}^1 (R,V^1)\Rho 1)\subset\gamma\subset
    \operatorname{emb}^1 (R,V^1)\Rho 1\}.
  \end{align*}

The inductive stage of the construction is this.  For $2\le m\le n$, given $\mathcal U^{m-1}\subset\otimes_{j=1}^n\mathcal S$,
we apply Lemma~\ref{l.few-small} to get a set   $V^m$ satisfying
\begin{equation*}
V^m\supset \sh{\mathcal U^{m-1}},\qquad \abs {V^m}\le{}(1+\delta){}
\abs{   \sh{\mathcal U^{m-1}}}.
\end{equation*}
The embedding function for rectangles $R\in\mathcal U^{m-1}$  is
 \begin{align*}
 \operatorname{emb}^m (R,V^{m}):=\sup\{\mu\ge1\mid &
 R_{1}\times\cdots\times{}R_{m-1}\times\mu R_{m}
 \\&\qquad\times\Rho {m+1}\times\cdots\times\Rho n\subset{}V^{m}\}.
 \end{align*}
 And the conclusion of (\ref{e.few-small}) holds.
 The collection $\mathcal U^m$ is then taken to consist of all rectangles of the form
 \begin{equation*}
 \otimes_{j=1}^{m-1} \Rho j\times \gamma\times \otimes_{j=m+2}^n\Rho j
 \end{equation*}
 where $R\in\mathcal U^{m-1}$ and $\gamma\in\mathcal S$  satisfies
 \begin{equation*}
( \Rho m\cup  \tfrac14\operatorname{emb}^m (R,V^{m}) \Rho m)\subset\gamma\subset \operatorname{emb}^m (R,V^{m})\Rho m.
 \end{equation*}

 To prove Lemma~\ref{l.journe}, we take $V:=V^n$.   It is the case that
 \begin{align*}
 \abs{V^{n}}\le{}&(1+\delta)\abs{\sh {\mathcal U^{n-1}}}
 \\{}\le{}&(1+\delta)\abs{V^{n-1}}
 \\{}\le{}&(1+\delta)^{n}\abs{\sh {\mathcal U}}.
\end{align*}

The definition of the embedding function is not so straight forward.  It is taken to be
\begin{equation*}
\operatorname{emb}(R)=\tfrac1{16}\inf_{1\le{}m\le{}n}\beta^m(R)
\end{equation*}
where $\beta^m(\cdot)$ are inductively defined below. The function $\imath(R)$ is taken to be the coordinate in which the infimum for the embedding function is achieved.

Set $\beta^1(R):=\operatorname{emb}^1 (R,V^{1})$.
In the inductive step,  for $2\le m\le n$, set $
\gamma_m(R):=\inf_{j<m}\beta^j(R)$.  For $1<\gamma<\gamma_m(R)$, let
\begin{equation*}
\beta^m_\gamma(R):=\operatorname{emb}^m (\varphi^m_\gamma(R),V^{m})
\end{equation*}
where $\varphi^m_\gamma(R)\in\mathcal U^{m-1}$ is the rectangle with $\varphi^m_\gamma(R)_{j}=\Rho j$ for $j\ge{}m$, and for
$1\le{}j<m$, $\varphi^m_\gamma(R)_{j}$ is the element of $\mathcal S$ of maximal length such that
\begin{equation*}
(\Rho j\cup\tfrac14\gamma\Rho j)\subset \varphi^m(R)_{j}\subset \gamma\Rho j.
\end{equation*}
Now,  take $\overline\gamma$ to be the largest value of $1\le{}\gamma\le{}\gamma_m(R)$ for which we have the inequality
$
\beta^m_\gamma(R)\ge\gamma
$.  Let us see that this definition of $\overline \gamma$ makes sense.
This last inequality is strict  for $\gamma=1$,  and as $\gamma$ increases, $\beta^m_\gamma(R)$ decreases, so $\overline\gamma$
is a well defined quantity.
Then define $\beta^m(R):=\beta^m_{\overline \gamma}(R)$, and for our use below,
set $\varphi^m(R):=\varphi^m_{\overline\gamma}(R)$.

The choices above prove our lemma, as we show now.
For each rectangle $R\in\mathcal U$, it is clear that $\operatorname{emb}(R)R\subset V$.
Take $\mathcal U'\subset \mathcal U$.
If we consider the sets $F(I,k,m,\mathcal U')$ as in (\ref{e.F}), then, by Lemma~\ref{l.few-small} applied in the $m$th coordinate,
\begin{equation*}
\sum_{I\in\mathcal D}\abs{F(I,k,m,\mathcal U')}\le{}
2^{\epsilon k}\abs{\sh {\varphi^m(\mathcal U')}}.
\end{equation*}
While we have a very good estimate for the shadow of $\varphi^m(\mathcal U)$, a corresponding good estimate for
an arbitrary subset $\mathcal U'$ seems very difficult to obtain.  But  it is a consequence of our construction
  that the rectangle $\varphi^m(R)$ is a rectangle which agrees with $R$ in the coordinates $j\ge{}m$ and, for
coordinates $1\le{}j<m$, is expanded by at most $32\operatorname{emb}(R)\le{}2^{k+6}$.  Hence, we have the estimate
\begin{equation*}
\ABs{\bigcup\{ \varphi^m(R)\mid R\in \mathcal U'\}}\lesssim{}2^{n k}\abs{\sh {\mathcal U'}}.
\end{equation*}
This follows from the weak $L^1$ bound for the maximal function in one dimension, applied in each coordinate separately.
It is in this last step that we lose the large power of the embeddedness.  Our proof is complete.


\subsection*{The $BMO_{-1} (\otimes_1^n\mathbb C_+)$ Formulation}

The form in which we apply the previous lemma is this.

\begin{lemma}\label{l.journe-n-1}
  Given a function $b$ with finite $BMO_{-1} (\otimes_1^n\mathbb C_+)$ norm and a collection of rectangles $\mathcal U$ whose
shadow has finite measure, the following construction is possible.  For all $\epsilon,\delta>0$, there is a set $V\supset\sh {\mathcal U}$
with $\abs V\le{}(1+\delta)\abs{\sh {\mathcal U}}$.  To each $R\in\mathcal U$, there is a quantity   $\Emb R..\ge1$ so that
\begin{gather*}
\Emb R..R\subset{}V, \qquad R\in\mathcal U,
\\
\NOrm \sum_{R\in\mathcal U} {\Emb R..}^{-(n+\epsilon)}{\ip b,v_R.v_R }.BMO (\otimes_1^n\mathbb C_+).\le{}K_{\delta,\epsilon}\norm b.BMO_{-1} (\otimes_1^n\mathbb C_+)..
\end{gather*}
\end{lemma}

For the proof, we apply Lemma~\ref{l.journe}.   To check the conclusion of the lemma, we take a subset $\mathcal U'\subset\mathcal U$
consisting of rectangles with $2^k\le{}\Emb R..<2^{k+1}$.  We then have
\begin{align*}
\sum_{R\in\mathcal U'} \abs{\ip b,v_R.}^2\le{}&\norm b. BMO_{-1} (\otimes_1^n\mathbb C_+).^2\sum_{I\in\mathcal D}\sum_{m=1}^n
\abs{ F(I,k,m,\mathcal U')}
\\\lesssim{}&\norm b. BMO_{-1} (\otimes_1^n\mathbb C_+).^2 2^{(n+\epsilon) k}\abs{\sh{\mathcal U'}}.
\end{align*}
This is all we need to  prove Lemma~\ref{l.journe-n-1}.

\end{document}